\documentclass[12pt]{amsart}
\usepackage{latexsym, amsmath,amssymb}

\usepackage{graphicx}

\usepackage{hyperref}

 \numberwithin{equation}{section}

\def\XXint#1#2#3{{\setbox0=\hbox{$#1{#2#3}{%
\int}$ }
\vcenter{\hbox{$#2#3$ }}\kern-.6\wd0}}

\renewcommand{\epsilon}{\varepsilon}
\newtheorem{theorem}{Theorem}

\newtheorem{lemma}[theorem]{Lemma}
\newtheorem{corr}[theorem]{Corollary}

\newtheorem{proposition}[theorem]{Proposition}
\newtheorem{deff}[theorem]{Definition}

\newcommand{\bth}{\begin{theorem}}
\newcommand{\ble}{\begin{lemma}}
\newcommand{\bcor}{\begin{corr}}

\newcommand{\bdeff}{\begin{deff}}

\newcommand{\bprop}{\begin{proposition}}
\newcommand{\ele}{\end{lemma}}
\newcommand{\ecor}{\end{corr}}
\newcommand{\edeff}{\end{deff}}

\numberwithin{theorem}{section}

\newcommand{\eprop}{\end{proposition}}

\renewcommand{\Pi}{\varPi}

\renewcommand{\epsilon}{\varepsilon}

\begin{document}

\title[C\'ordoba's differentation theorem: revisited]
{C\'ordoba's differentiation theorem: revisited}

\author[A. D. Mart\'inez]{\'Angel D. Mart\'inez}
\address{Institute for Advanced Study, Fuld Hall 412, 1 Einstein Drive, Princeton, NJ 08540, United States of America} 
\email{amartinez@ias.edu}


\begin{abstract}
In this paper we prove an exponential covering lemma implying the three dimensional case of a well-known conjecture formulated by A. Zygmund circa 1935 and solved by A. C\'ordoba in 1978. Our approach avoids a subtle argument involving the power series of the exponential function. 
\end{abstract}

\maketitle
\section{Introduction}

The concept of maximal function has been  central in mathematical analysis since its appeareance in Hardy  and Littlewood's seminal work \cite{HL}. It allows to obtain the so-called Lebesgue's differentiation theorem in a quantitative way. In 1935 Jessen, Marcinkiewicz and Zygmund proved that instead of balls or cubes one may also differentiate using the class $\mathcal{B}_n$, of shrinking $n$-dimensional rectangles (paralellepipeds) with sides parallel to the axes, provided the function is locally in the class $L^1(1+\log_+L^1)^{n-1}(\mathbb{R}^n)$. The quantitative version of this result is known as the strong maximal theorem: let us define
\[M_nf(x)=\sup_{x\in R\in\mathcal{B}_n}\frac{1}{m(R)}\int_R|f(y)|dm(y)\]
where $f$ is a locally integrable function and $m$ denotes Lebesgues's measure in $\mathbb{R}^n$. Then we have:
\begin{theorem}[Jessen, Marcinkiewicz and Zygmund, \cite{JMZ}]\label{mst}
The class $\mathcal{B}_n$ differentiates functions $f$ which are locally in $L^1(1+\log_+L^1)^{n-1}(\mathbb{R}^n)$, that is, we have
\[\lim_{\textrm{diam}(R)\rightarrow 0}\frac{1}{m(R)}\int_Rf(y)dm(y)=f(x)\textrm{ for a.e. $x\in\mathbb{R}^n$}\]
where the limit is taken over all rectangles such that $x\in R\in\mathcal{B}_n$. Quantitatively, the following inequality
\[m\{x\in\mathbb{R}^n:M_nf(x)>\alpha\}\leq C_n\int_{\mathbb{R}^n}\frac{f(x)}{\alpha}\left(1+\log_+\frac{f(x)}{\alpha}\right)^{n-1}dm(x)\]
holds for some constant $C_n$ that depends only on the dimension.
\end{theorem}
Their proof relies on an iteration scheme based on the one dimensional result and making use of the product structure of $\mathbb{R}^n$. Quite surprisingly the Orlicz space, $L^1(1+\log_+L^1)^{n-1}(\mathbb{R}^n)$, is sharp if one assumes all its functions to be differentiated by $\mathcal{B}_n$ (cf. \cite{JMZ}, theorem 8). It should be pointed out that the quantitative part of the statement appeared later (cf. \cite{Fava} or \cite{guz}).

Nowadays the standard proof of the weak $(1,1)$ inequality for the Hardy-Littlewood maximal function (i.e. the case of cubes or balls in $\mathbb{R}^n$) relies on the Vitali covering lemma which allows to deduce this estimates as a consequence of the geometry of balls. Standard references are the treatises of de Guzm\'an \cite{G} and Stein \cite{St} where $C_n=5^n$. 

It was not until 1975 when the exponential type covering needed to prove geometrically the strong maximal theorem was finally discovered (cf. \cite{ACF}). The proof obtained by the A. C\'ordoba and R. Fefferman makes use of certain sparseness properties. Let us introduce the one we shall need, using their notation: 
\begin{itemize}
\item[($P_1$)] A sequence of rectangles $\{R_j\}$ satisfies the $P_1$ sparseness property if for any $k$ the following inequality holds: 
\[m\left(R_k\cap\bigcup_{j<k}R_j\right)\leq\frac{1}{2}m(R_k).\]
\end{itemize}

In 1978 A. C\'ordoba took advantage of the geometric approach to settle Zygmund's conjecture in the three-dimensional case:

\begin{theorem}[A. C\'ordoba, \cite{AC2, AC3, AC4}]\label{acor}
Given any function $\phi$ non decreasing in each variable separately let $\mathcal{B}_{\phi}(\mathbb{R}^3)$ consist of rectangles whose lengths are of the form $(x,y,\phi(x,y))$. Then it differentiates $L^1(1+\log_+L^1)(\mathbb{R}^3)$.
\end{theorem}

In fact, C\'ordoba's differentiation theorem is stronger but let us state it in this way for the sake of the exposition's clarity. His elegant proof relies heavily on the three dimensional setting. No purely analytical proof has been found to-date. For applications of this result we refer the reader to the original works of A. C\'ordoba. Zygmund's curiosity was not satisfied by this and asked for further generalizations of this result in arbitrary dimension. Some of them were dismissed in \cite{So} by F. Soria. His counterexamples have side lenghts of  the form $(s,\phi(s,t),\psi(s,t))$ where the functions $\phi$ and $\psi$ are monotone increasing. Recently more counterexamples have been found by G. Rey in \cite{Rey}. The rest remain an intriguing open question. 

In the next section we provide a geometric covering result that implies Theorem \ref{acor}, we refer the reader to the literature for this implication (cf. \cite{AC2, AC3, ACF}). Our approach neatly shows why and where the argument breaks down in the higher dimensional case.

\section{Proof of Theorem \ref{acor}}\label{monotone}

Let us now state the main result of this paper:

\begin{theorem}[A. C\'ordoba, 1978]
Given any function $\phi$ non decreasing in each variable separately let $\mathcal{B}_{\phi}(\mathbb{R}^3)$ consist of rectangles whose lengths are of the form $(x,y,\phi(x,y))$. There exist dimensional constants $C, c>0$, but otherwise independent, such that for any $\{R_{\alpha}\}\subseteq\mathcal{B}_{\phi}(\mathbb{R}^3)$ there exists a subsequence $\{R_j\}$ such that
\begin{itemize}
\item[(1)] 
\[m\left(\bigcup_{\alpha}R_{\alpha}\right)\leq C m\left(\bigcup_{j}R_j\right).\]
\item[(2)]
\[\int_{\bigcup_jR_j}\exp c \left(\sum_j\chi_{R_j}\right)^{}dm\leq C m\left(\bigcup_jR_j\right).\]
\end{itemize}
\end{theorem}

\textsc{Proof:} it is well known that without loss of generality we may assume that {\em a priori} the rectangles $\{R_{j}\}$ to be a finite sequence ordered by decreasing length of their third side length and satisfying the $P_1$ property. This reduction respects (1) from the statement (cf. \cite{ACF}). We will denote by $R^*$ a rectangle centered in the same point and dilated by a factor of three. The introduction of this will be apparent later when we require intersections to be clean.

We select the rectangles following the following algorithm. Having chosen $k-1$ of them, we will choose $R_k$ to be the next enlisted rectangle $R$ that satisfies
\begin{equation}\label{localmean}
\frac{1}{m(R)}\int_{R}\exp \left(\sum_{j=1}^{k-1}\chi_{R^*_j}\right)^{}dm\leq 3.
\end{equation}
This can be shown to imply property (2) from the statement. Denoting
\[I_k=\int_{\bigcup_{j=1}^kR_j}\textrm{exp }\left(\sum_{j=1}^k\chi_{R^*_j}\right)^{}dm\]
let us observe that
\[I_k=\int_{\bigcup_{j=1}^{k-1}R_j\setminus R_k}\textrm{exp }\left(\sum_{j=1}^{k-1}\chi_{R^*_j}\right)^{}dm+\int_{R_k}\textrm{exp }\left(\sum_{j=1}^k\chi_{R^*_j}\right)dm.\]
The choice of of $R_k$ shows that $I_k\leq I_{k-1}+3e m(R_k)$ which  by induction implies $I_k\leq 3e\sum_{j=1}^km(R_j)$. From this and the $P_1$ property our claim follows easily, namely
\begin{equation}\label{3}
\int_{\bigcup_{j=1}^kR_j}\textrm{exp} \left(\sum_{j=1}^k\chi_{R^*_j}\right)^{}dm\leq 6 e m\left(\bigcup_jR_j\right).
\end{equation}
This shows that, following this sieve, part (2) will hold with $c=1$  and it only rests to show that it also respects (1). To do so we observe that for a rejected rectangle $R$ the negation of the inequality \ref{localmean} implies
\[3<\sum_{k=0}^{\infty}A_ke^k\]
holds, where we are using the notation $A_0=1$ and
\[A_k=\frac{m(\{x\in R:\sum_j\chi_{R^*_j}(x)=k\})}{m(R)}\]
for $k>0$. Notice that $A_k$ is the proportion of points in $R$ that belong to exactly $k$ rectangles. Let us denote by $R^1_j$ those rectangles previous to $R$ in the selection process whose first and third side lengths are greater to those of $R$. Analogously, from the remaining rectangles, we denote by $R^2_j$ those with second and third lengths exceeding those of $R$. Decomposing $A_k$ in different cases (namely, the $k=r+s$ intersections being the result of the intersection of $r$ rectangles from the first class and $s$ from the second) the above inequality is equivalent to
\[3<\sum_{k=0}^{\infty}\sum_{r+s=k}A_{r,s}e^k\]
where 
\[A_{r,s}=\frac{m(\{x\in R:\sum_j\chi_{R^{1*}_j}(x)=r\textrm{ and }\sum_j\chi_{R^{2*}_j}(x)=s\})}{m(R)}.\]
\begin{figure}[htb]
\centering
\includegraphics[width=90mm]{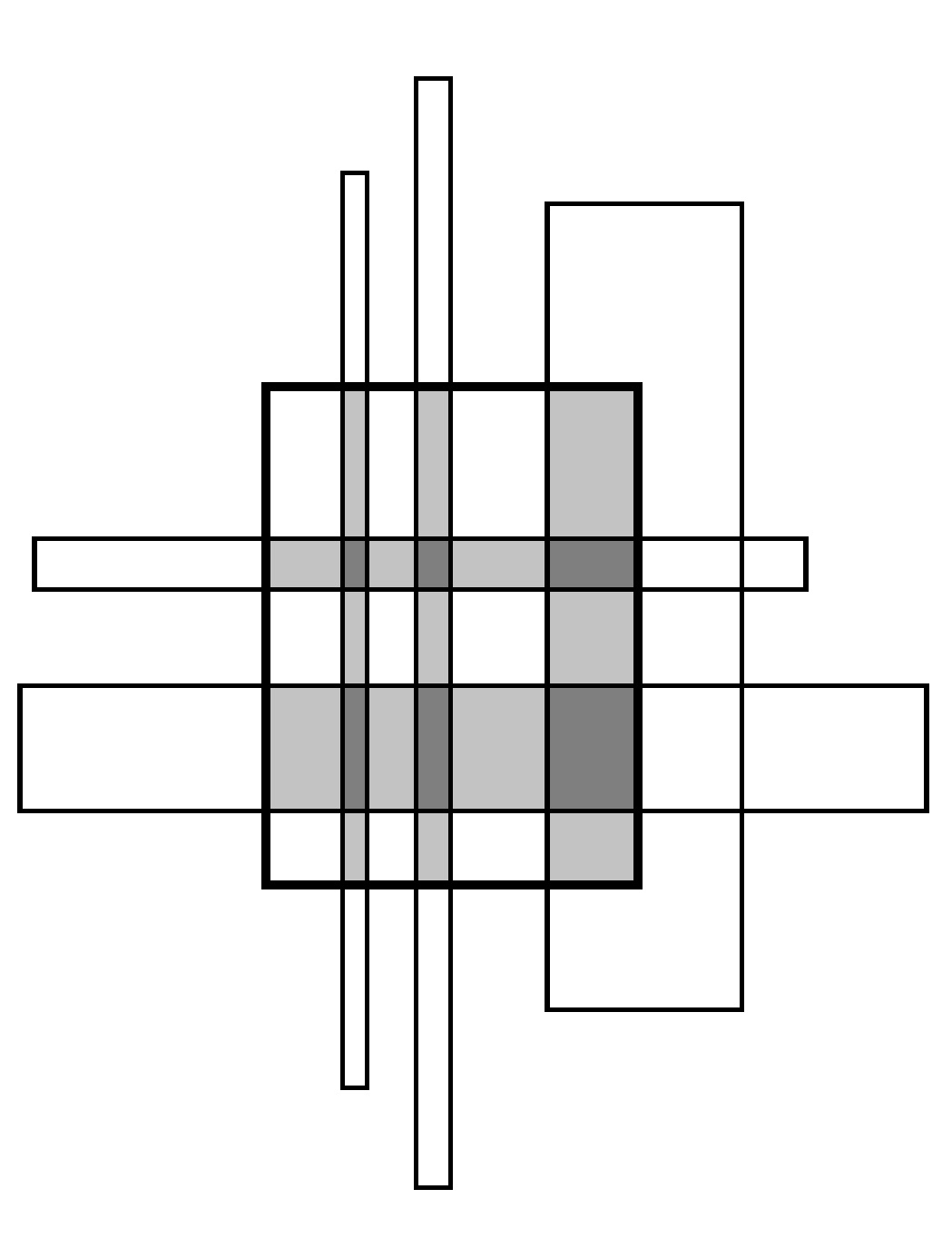}
\caption{Schematic $(x_1,x_2)$-section of the intersection of $R$ with three vertical and two horizontal rectangles (respectively $R^{1*}_j$ and $R^{2*}_j$). In different shades of gray $A_{1,0}$, $A_{0,1}$ (lighter) and $A_{1,1}$ (darker).}\label{fig:b1}
\end{figure}
Let us introduce a related piece of notation now. Let $a_0=1$ and in general
\[a_r=\frac{m(\{x\in R:\sum_j\chi_{R^{1*}_j}(x)=r\})}{m(R)}\]
for $r>0$. Similarly, let $b_0=1$ and
\[b_s=\frac{m(\{x\in R:\sum_j\chi_{R^{2*}_j}(x)=s\})}{m(R)}\]
for $s>0$. It is now obvious from the clean intersections of $R$ and $R^{\ell*}_j$, $\ell=1,2$, (see Fig. 1) that $A_{r,s}\leq a_rb_s$ holds from which it follows that
\[3<\left(\sum_{r=0}^{\infty}a_re^r\right)\left(\sum_{s=0}^{\infty}b_se^s\right).\]
(At this point the proof breaks down in higher dimensions.) Reading back what this means it implies
\[\max_{\ell=1,2}\left\{\frac{1}{m_{x_{\ell}}(R)}\int_R \exp\left(\sum_{i=1}^k\chi_{R^*_i}\right)dm_{x_{\ell}}\right\}>\sqrt{3}-1\]
where $m_{x_\ell}$ is the length in the $\ell$th direction. The integral mean is a maximal function of the exponential restrected to any line intersecting $R$ and parallel to the $x_{\ell}$ axis. We will denote such a (one dimensional) maximal operator by $M_{\ell}$. The proof concludes estimating the size of the set covered by rejected rectangles. For this purpose one observes that it is included in
\[\left\{x: M_{1}\exp\left(\sum_{i=1}^k\chi_{R^*_i}\right)(x)>\sqrt{3}-1\textrm{ or }M_{2}\exp\left(\sum_{i=1}^k\chi_{R^*_i}\right)(x)>\sqrt{3}-1\right\}\]
whose measure can be bounded appropiately by the weak $L^1$ boundedness of the Hardy-Littlewood maximal function. Indeed, the maximal functions appearing above are one dimensional means in the first and second coordinates alone to which the weak $L^1$ bound applies linewise. Integrating this $L^1$ bounds on the rest of variables and using Theorem \ref{mst} for $n=1$ one gets
\[m\left(\bigcup_{\beta}R_{\beta}\right)\leq \frac{5}{\sqrt{3}-1}\int_{\bigcup_{R_i}}\exp\left(\sum_{i=1}^k \chi_{R^*_i}\right)dm\]
where $\beta$ extends over the indices of rejected rectangles. This can be bounded using equation \ref{3} to conclude (1).

\textsc{Remark:} it is easy to show that this argument also works for a family satisfying that in some order any given $R$ has two sidelengths smaller than any previous rectangle.

\section{Acknowledgments}

The author would like to express his gratitude to Antonio C\'ordoba for introducing him to the subject in a graduate course delivered in Fall 2014 at UAM-ICMAT and Eric Latorre for his (unpublished) lecture notes \cite{ACE}. The author is indebted to the anonymous referees for their careful reading of the manuscript, suggestions and for pointing out references \cite{Fava} and \cite{guz} to his attention.

The author was partially supported by  MTM2014-56350-P  project of the MCINN (Spain). This material is based upon work supported by the National
Science Foundation under Grant No. DMS-1638352. 

\nocite{*}

\section{Acknowledgments}

The author would like to express his gratitude to Antonio C\'ordoba for introducing him to the subject in a graduate course delivered in Fall 2014 at UAM-ICMAT and Eric Latorre for his (unpublished) lecture notes \cite{ACE}. The author is indebted to the anonymous referees for their careful reading of the manuscript, suggestions and for pointing out references \cite{Fava} and \cite{guz} to his attention.

The author was partially supported by  MTM2014-56350-P  project of the MCINN (Spain). This material is based upon work supported by the National
Science Foundation under Grant No. DMS-1638352.

\end{document}